\documentclass{article}
\usepackage[british]{babel}
\usepackage{amsmath,amssymb,enumerate,latexsym,theorem}
\usepackage[parfill]{parskip}
\usepackage{longtable}
\usepackage{cite}

\newcommand{\mat}[0]{\operatorname{Mat}}
\newcommand{\herm}[0]{\operatorname{Herm}}
\newcommand{\sq}[0]{\operatorname{Sq}}

\newcommand{\matbar}[0]{\overline\mat}
\newcommand{\hermbar}[0]{\overline\herm}
\newcommand{\sqbar}[0]{\overline\sq}

{\theorembodyfont{\rmfamily}}

\begin{document}

\title{Towards a Singular Value Decomposition and spectral theory for all rings}
\author{Ran Gutin\footnote{Department of Computer Science, Imperial College London}}
\date{\today}

\maketitle

\begin{abstract}
  We propose definitions of SVD, spectral decomposition (for self-adjoint matrices) and Jordan decomposition which make sense for all rings. For many rings, these decompositions can be shown to exist. For some specific rings, these decompositions are complicated to describe in full and prove the existence of. These decompositions have occurred piecemeal in the literature. We conjecture that they exist for many rings, including all Clifford algebras over the real numbers and complex numbers. The origin of this programme is not directly in module theory or linear algebra.
\end{abstract}

\section{Motivation before giving the
definitions}\label{motivation-before-giving-the-definitions}

\subsection{Objective}\label{objective}

The objective of this paper is to propose a general definition for three
matrix decompositions, which can be shown to be satisfiable (or not)
over any given ring:

\begin{itemize}
\item
  SVD (Singular Value Decomposition)
\item
  Spectral decomposition. In this paper, we see this as a decomposition
  of self-adjoint matrices.
\item
  Jordan decomposition.
\end{itemize}

The motivation for this is that we have discovered analogues of those
three decompositions outside the obvious linear algebra or module theory
context. These analogues present canonical forms which are not as simple
as in the complex case, because they might fail to be diagonal. This
phenomenon is already familiar for the Jordan decomposition (in the
complex case), but is not present for the other two decompositions when
working over the complex numbers. In general, the canonical forms can be
complicated, and their existence is not easy to prove on a case-by-case
basis.

\subsection{Summary of origins}\label{summary-of-origins}

I give an account of how I stumbled upon SVDs over different rings.
Admittedly, this is somewhat subjective. It is written in the first
person.

There is a 19th century ``non-Euclidean'' geometry (in the Kleinian
sense) called \emph{Laguerre geometry} which admits a relationship to
the algebra of \(2 \times 2\) matrices over the dual numbers. I found
that based on the correspondence between the \emph{congruence
transformations} of this geometry, and the \(2 \times 2\) matrices over
the dual numbers, there was a strong hint towards the existence of an
analogue of the Singular Value Decomposition over the ring of dual
numbers instead of the usual ring of complex numbers. An account of this
matrix/transformation correspondence can be found in Yaglom's book
\emph{Complex numbers in geometry}. The SVD interpretation and
subsequent result was inspired by a classification of the Laguerre
transformations which can be found in Yaglom's book. This analogue of
the SVD turned out to exist for matrices over the dual numbers, where
the matrices could be of any possible dimension (not just the
\(2 \times 2\) case that Yaglom considered, and not just the square
matrices), even in the case where the matrices were singular, and
satisfied certain uniqueness properties similar to the complex SVD.

I will now state exactly the result I obtained \cite{dualsvd}. The ring
of dual numbers is denoted formally as
\(\mathbb R[\varepsilon]/(\varepsilon^2)\).
\footnote{More informally, a dual number is a number of the form $a + b \varepsilon$ where $\varepsilon^2 = 0$ while $\varepsilon \neq 0$. The dual numbers form an associative and commutative unital algebra over the real numbers. }The
dual numbers admit an involution \(*\) for which
\((a + b\varepsilon)^* = a - b \varepsilon\). This defines a
conjugate-transpose or adjoint operator on matrices \(M^*\) such that
\((M^*)_{ij} = (M_{ji})^*\). Analogously, we get a notion of unitary
matrix \(U\) for which \(UU^* = U^* U = I\). Given a matrix \(M\) over
the dual numbers (of any dimension), I obtained a result which states
that \(M\) can be factorised into \(USV^*\) where \(U\) and \(V\) are
unitary, and \(S\) is a direct sum of matrices in the set

\[
  \begin{aligned}
  &\{[x] : x \in \mathbb R; x > 0\}\\& \cup \left\{\begin{bmatrix} x & -y\varepsilon \\ y\varepsilon & x\end{bmatrix} : x \in \mathbb R, y \in \mathbb R; x > 0, y > 0\right\}\\& \cup \{[y\varepsilon] : y \in \mathbb R; y > 0\}\\&\cup \{0_{1,0}, 0_{0,1}\}.
  \end{aligned}
\]

This set is quite complicated, and we won't dwell on it. Notice though
that the matrix \(0_{1,0}\) (and \(0_{0,1}\)) is missing in the usual
account of the SVD over \(\mathbb C\), but we will show that it is implicitly
there as well!
It's the unique \(1 \times 0\) matrix.

To explain the occurrence of this strange matrix, let's consider the
complex matrix
\(M = \begin{bmatrix}1 & 2 & 0 \\ 2 & 1 & 0 \end{bmatrix}\). Ignoring
the exact values of \(U\) and \(V\), but focussing our attention on
\(S\), we have that \(M = USV^*\) where
\(S = [3] \oplus [-1] \oplus 0_{0,1}\). The reader should verify that
given a matrix \(K\), the value of \(K \oplus 0_{0,1}\) is the same
matrix as \(K\) but padded with a zero column. Likewise, the value of
\(K \oplus 0_{1,0}\) is the same matrix as \(K\) but padded with a zero
row. The value of \([3] \oplus [-1] \oplus 0_{0,1}\) is thus
\(\begin{bmatrix}3 & 0 & 0 \\ 0 & -1 & 0 \end{bmatrix}\) as expected. We
conclude that the ``singular values'' of a complex matrix form a multiset
(a similar object to a \emph{set} but the elements are allowed to
repeat) whose elements belong to the set
\(\{[x] : x \in \mathbb R; x > 0\} \cup \{0_{1,0}, 0_{0,1}\}\). This is
actually a set of \emph{matrices} rather than scalar ``values''. This
suggests that the notion of a singular value as a scalar value is
questionable.

Given a matrix \(M\) over \(\mathbb C\) or
\(\mathbb R[\varepsilon]/(\varepsilon^2)\), the direct summands that
make up \(S\) in \(M = USV^*\) are unique up to permutation. In the dual
number case, a direct summand can be of dimension \(2 \times 2\), and
non-square direct summands can be found over both rings (and in fact all
rings).

Upon discovering the above dual number SVD, I found that a dual number
SVD of sorts had already been considered in the literature, and it
wasn't the one above. In that case, the involution over the dual numbers
is the trivial one: \(z^* := z\), and the notion of unitary matrix
degenerates into \(U^T U = UU^T = I\). We then have that the direct
summands that make up \(S\) are of the form:

\[\{[x + y\varepsilon] : x \in \mathbb R, y \in \mathbb R; x > 0\} \cup \{[y \varepsilon] : y \in \mathbb R; y > 0\}\cup \{0_{1,0}, 0_{0,1}\}.\]

This shows that the correct setting is not a \emph{ring} as such, but a
\emph{\(*\)-ring}, which is a ring equipped with an involution. The
motivation for considering matrix algebra over this particular
\(*\)-ring came from the mechanics literature where this SVD was
suggested for solving kinematic synthesis problems.

Over the \(*\)-ring of \emph{double numbers}, which are defined as
\(\mathbb R \oplus \mathbb R\), equipped with the involution
\((a,b)^* := (b,a)\), a matrix decomposition that constitutes an SVD
over the double numbers has been given \cite{contragredient}. I attempted to come up with such
a decomposition myself, which I called the Jordan SVD, but it had the
problem that it didn't exist for all matrices over the double numbers,
and therefore didn't satisfy the general criteria I give here. The paper
in which the \emph{actual} double-number SVD is introduced has some
limitations, chief among them is that the paper is unaware of the
connection to the double numbers, or that the decomposition constitutes
\emph{the} analogue of the SVD over the \(*\)-ring of double numbers. A
lot of terminology is introduced there, like \emph{contragredient
equivalence}, which I view as needless.

Finally, over the \emph{quaternions} an analogue of SVD exists which
satisfies my definition. This has recently been extended to the case of
the \emph{dual quaternions} (an algebra equal to
\(\mathbb H[\varepsilon]/(\varepsilon^2)\)). These algebras are
non-commutative, so our definitions don't assume commutativity.

There is clearly a need for a general result which states precisely when
an analogue of the SVD, spectral decomposition for self-adjoint
matrices, or Jordan decomposition, exists for a given \(*\)-ring. Such a
result has already been obtained for the analogue of the Jordan
decomposition which we discuss here, but this appears at the moment to
be the least general of the three decompositions. We supply general
definitions here and formulate conjectures.

\subsection{Rough desired criteria}\label{rough-desired-criteria}

If \(M = USV^*\) is the ``SVD'' of M, then \(S\) should be a block
diagonal matrix which is unique up to permutation of the blocks. This is
similar to, and somewhat related to, the fact that the Jordan
decomposition of a square matrix \(M\), expressed as \(PJP^{-1}\),
should be unique up to permutation of the Jordan blocks. The objective
is to decide whether or not a set of such blocks can be exhibited for a
given ring.

The formal definitions which we give below express in a rigorous way the
desired existence and uniqueness properties which we've alluded to.

\section{Definitions}\label{definitions}

A \(*\)-ring is a ring \(R\) equipped with an automorphism denoted
\(*:R \to R\) that has order \(2\). Such an automorphism is called an
\emph{involution}. Every ring can be made into a \(*\)-ring in at least
one way by defining \(z^* = z\) for all \(z \in R\).

We define the following monoids:

\begin{itemize}
\item
  \(\mat(R,*)\) is the monoid of all matrices over a \(*\)-ring \(R\) where
  the monoid operation is \(\oplus\), denoting direct sum of matrices.
  The matrices can be of any possible dimensions, and they don't have to
  be square. The matrices can also have \(0\) rows or \(0\) columns.
\item
  \(\herm(R,*)\) is the monoid of all self-adjoint matrices over a
  \(*\)-ring \(R\) where the monoid operation is \(\oplus\).
\item
  \(\sq(R)\) is the monoid of all square matrices over a ring \(R\), where
  \(R\) is merely a ring and not a \(*\)-ring.
\end{itemize}

We then define the following equivalence relations on \(\mat(R,*)\) and its
various submonoids above:

\begin{itemize}
\item
  \({\sim_{\text{UE}}}\) means \emph{unitary equivalence}. In other
  words, we have that \(A {\sim_{\text{UE}}} B\) is true whenever there
  exist unitary matrices \(U\) and \(V\) such that \(A = UBV^*\).
\item
  \({\sim_{\text{US}}}\) means \emph{unitary similarity}. In other
  words, we have that \(A {\sim_{\text{US}}} B\) is true whenever there
  exists a unitary matrix \(V\) such that \(A = VBV^*\).
\item
  \({\sim_{\text{S}}}\) means \emph{similarity}. In other words, we have
  that \(A {\sim_{\text{S}}} B\) is true whenever there an invertible
  matrix \(P\) such that \(A = PBP^{-1}\).
\end{itemize}

We aim to study the three monoids:

\begin{itemize}
\item
  \(\matbar(R,*) := \mat(R,*) / {\sim_{\text{UE}}}\) with the intention of generalising the
  \emph{singular value decomposition},
\item
  \(\hermbar(R,*) := \herm(R,*) / {\sim_{\text{US}}}\) with the intention of generalising the
  \emph{spectral theorem} (on self-adjoint matrices).
\item
  \(\sqbar(R,*) := \sq(R) / {\sim_{\text{S}}}\) with the intention of generalising the
  \emph{Jordan decomposition}.
\end{itemize}

To do this, notice that all three monoids are abelian, and for some
\(*\)-rings \(R\), they are even \emph{free} as abelian monoids. A
\emph{free abelian monoid} consists of all finite multisets whose
elements belong to some set \(S\). An isomorphism between each of the
three monoids above and a free abelian monoid produces an analogue of
the SVD, spectral theorem, and Jordan decompositions respectively.

As an aside: Note that all six parametrised families of monoids can be
thought of as \emph{functors} between two categories, if this fact can
ever be useful. The functors are all from the category of \(*\)-rings to
the category of monoids.

\section{Discussion of the relationship of monoid algebra to the SVD,
spectral decomposition and Jordan
decompositions}\label{discussion-of-the-relationship-of-monoid-algebra-to-the-svd-spectral-decomposition-and-jordan-decompositions}

When an abelian monoid \((M,+,0)\) is free, there is a \emph{unique}
subset (which we will here denote \(P\)) of \(M\) (called the generators
of \(M\)) such that each element of \(M\) is a unique sum elements of
the generators \(P\). It can be argued that many ``unique
factorisation'' type results in mathematics merely state the fact that
some abelian monoid with a complicated construction is actually free.

An example of the above is the abelian monoid of \emph{positive
integers} under integer multiplication, in which the generators are the
prime numbers. The fact that for every integer \emph{there exists} a
prime factorisation, and moreover it is \emph{unique}, is equivalent to
the fact that the abelian monoid of positive integers is free.

The abelian monoid \(\matbar(R,*)\) for a given
\(*\)-ring \((R,*)\) is constructed in a complicated way. In spite of
that fact, if it is free, then it is actually quite simple. Its freeness
captures the existence and uniqueness of the SVD.

\section{Summary of existing
results}\label{summary-of-existing-results}

\textbf{Theorem 1}: \(\sqbar(R)\) is a free abelian
monoid for every Artinian ring \(R\).

\emph{Proof}. See \cite{ringjordan}.

The above theorem in particular shows that something approaching a
Jordan decomposition exists for the dual numbers, which are an Artinian
ring. All finite dimensional associative algebras over a field, like the
quaternions and dual numbers, are Artinian, so this generalises an
important aspect of the Jordan decomposition.

Let \(\operatorname{swap}(a,b) := (b,a)\). This is a natural choice of
involution for \(R \oplus R\) where \(R\) is any ring. An illustrative
special case is when \(R\) is \(\mathbb R\), in which case
\(R \oplus R\) is commonly called either the \emph{split-complex
numbers} or the \emph{double numbers}. This is an interesting
hypercomplex number system considered in \cite{doublenumbers}. When
\(R = \mathbb C\), the name we give to \(R \oplus R\) is the
\emph{double complex numbers}.

\textbf{Theorem 2}:
\(\hermbar(R\oplus R,\operatorname{swap})\cong \sqbar(R)\)
for every ring \(R\).

\textbf{Remark}: When both sides are treated as functors of \(R\), the
isomorphism is natural. We don't attempt to verify this fact.

\emph{Proof}. Let \(M\) be an element of
\(\herm(R \oplus R, \operatorname{swap})\). We have that \(M = (A,B^T)\)
because \(M\) in particular must belong to \(J(R \oplus R)\). We
furthermore have that \(M\) is Hermitian, so
\((A,B^T) = M = M^* = (B,A^T)\); therefore \(A = B\). In summary, we
have that \(M = (A,A^T)\). Let
\(\phi_R : \sq(R) \,\to \herm(R \oplus R, \operatorname{swap})\) be given by
\(\phi_R(A) = (A,A^T)\). The mapping is clearly an isomorphism of
monoids. But we're not done yet because we would like an isomorphism
\(\psi_R : \sqbar(R) \to \hermbar(R \oplus R, \operatorname{swap})\)
instead. This isomorphism is obtained by noticing that if
\(M \sim_\text{S} K\) then \(\phi_R(M) \sim_\text{US} \phi_R(K)\). We
show that this is indeed true: If \(M {\sim_{\text{S}}} K\) then
\(M = PKP^{-1}\), so
\(\phi_R(M) = (M,M^T) = (PKP^{-1},(P^T)^{-1} K^T P^T) = (P,(P^{-1} )^T) (K,K^T) (P^{-1},P^T) {\sim_{\text{US}}} (K,K^T) = \phi_R(K)\).
We are done.

The above result is significant because it shows that the monoid family
\(\sqbar(R)\) is somewhat redundant. The study of these
monoids is subsumed by the study of monoids in the family
\(\hermbar(R,*)\).

\textbf{Theorem 3}: \(\matbar(R,*)\) admits an
isomorphism to a free abelian monoid with the following generators when
\((R,*)\) is any of the following \(*\)-rings:

  \begin{longtable}{|l|l|l|l|}
    \hline
  & \(R\) &  \(*:R \to R\) & Generators of \(\matbar(R,*)\)\\
  \hline
  \hline
  1 & \(0\) & \(z^* := z\) & \(\{0_{1,0}, 0_{0,1}\}\) \\
  \hline
  2 & \(\mathbb R\) & \(z^* := z\) &
  \(\begin{aligned}&\{[x] : x \in \mathbb R; x > 0\}\\ &\cup \{0_{1,0}, 0_{0,1}\}\end{aligned}\) \\
  \hline \\
  3 & \(\mathbb C\) & \((a + bi)^* := a - bi\) &
  \(\begin{aligned}&\{[x] : x \in \mathbb R; x > 0\} \\&\cup \{0_{1,0}, 0_{0,1}\}\end{aligned}\) \\
  \hline \\
  4 & \(\mathbb R[\varepsilon] / (\varepsilon^2)\) & \(z^* := z\) &
  \(\begin{aligned}&\{[x + y\varepsilon] : x \in \mathbb R, y \in \mathbb R; x > 0\}\\ &\cup \{[y \varepsilon] : y \in \mathbb R; y > 0\}\cup \{0_{1,0}, 0_{0,1}\}\end{aligned}\) \\
  \hline \\
  5 & \(\mathbb R[\varepsilon] / (\varepsilon^2)\) &
  \((a + b\varepsilon)^* := a - b\varepsilon\) &
  \(\begin{aligned}&\{[x] : x \in \mathbb R; x > 0\}\\ &\cup \left\{\begin{bmatrix} x & -y\varepsilon \\ y\varepsilon & x\end{bmatrix} : x \in \mathbb R, y \in \mathbb R; x > 0, y > 0\right\}\\ &\cup \{[y\varepsilon] : y \in \mathbb R; y > 0\}\cup \{0_{1,0}, 0_{0,1}\}\end{aligned}\) \\
  \hline \\
  6 & \(\mathbb C \oplus \mathbb C\) & \(\operatorname{swap}\) &
  \(\begin{aligned}&\{(J_m(re^{i\theta}),J_m(re^{i\theta})^T) : r > 0, \theta \in [0,\pi)\} \\&\cup \left\{\left(\begin{pmatrix}I_m & 0\end{pmatrix},\begin{pmatrix}0\\ I_m\end{pmatrix}^T\right) : m \in \mathbb N\right\}\\ &\cup \left\{\left(\begin{pmatrix}0\\ I_m\end{pmatrix},\begin{pmatrix}I_m & 0\end{pmatrix}^T\right) : m \in \mathbb N\right\}\\&\cup \{(J_m(0), I_m^T) : m \in \mathbb N\}\\& \cup \{(I_m, J_m(0)^T) : m \in \mathbb N\}\cup\{0_{1,0}, 0_{0,1}\}\end{aligned}\) \\
  \hline
  7 & \(\mathbb H\) & \(\begin{aligned}&(a + bi + cj + dk)^* :=\\& a - bi - cj - dk \end{aligned}\) &
  \(\{[x] : x \in \mathbb R; x > 0\} \cup \{0_{1,0}, 0_{0,1}\}\) \\
  \hline
  8 &
  \(\begin{aligned}&\operatorname{span}(1, i, \varepsilon j, \varepsilon k)\\&\subset \mathbb H[\varepsilon]/(\varepsilon^2)\end{aligned}\)
  &
  \(\begin{aligned}&(a + bi + c\varepsilon j + d \varepsilon k)^* :=\\& a - bi - c\varepsilon j - d \varepsilon k\end{aligned}\)
  &
  \(\begin{aligned}&\{[x] : x \in \mathbb R; x > 0\}\\& \cup \left\{\begin{bmatrix} x & -\delta \\ \delta & x\end{bmatrix} : x \in \mathbb R; x > 0, \delta^2 = 0\right\}\\& \cup \{[y\varepsilon j]: y \in \mathbb R; y > 0\}\cup \{0_{1,0}, 0_{0,1}\}\end{aligned}\) \\
  \hline
\end{longtable}
\emph{Proof}. For each example in turn:

\begin{enumerate}
\item
  The ring here is the \(0\) ring, whose set of elements is \(\{0\}\).
  The definition of the operations \(\{+,-,\times\}\) is immediate. The
  generators are the two degenerate matrices \(0_{1,0}\) and
  \(0_{0,1}\). The first matrix has \(1\) row and \(0\) columns. The
  second matrix has \(0\) rows and \(1\) column. Matrices with \(0\)
  rows or columns appear strange but they are an inevitable part of
  matrix algebra because they represent linear maps which map to or from
  \(0\)-dimensional vector spaces.
\item
  A proof of this fact can be found in any undergraduate textbook on
  linear algebra for mathematics students, as long as that textbook
  covers the spectral theorem and SVD. The only unusual generators are
  \(0_{1,0}\) and \(0_{0,1}\), which are \(1 \times 0\) and
  \(0 \times 1\) matrices respectively. While it is strange to have
  matrices with these dimensions, they are an inevitable part of the
  matrix formalism because they represent linear maps that go to or from
  \(0\)-dimensional vector spaces.
\item
  A proof of this fact can be found in any undergraduate textbook on
  linear algebra for mathematics students, as long as that textbook
  covers the spectral theorem and SVD. The only unusual generators are
  \(0_{1,0}\) and \(0_{0,1}\), which are \(1 \times 0\) and
  \(0 \times 1\) matrices respectively. While it is strange to have
  matrices with these dimensions, they are an inevitable part of the
  matrix formalism because they represent linear maps that go to or from
  \(0\)-dimensional vector spaces.
\item
  See \cite{dualsvd}. The only thing missing from
  the paper is a complete proof of uniqueness. This can be achieved by
  using the uniqueness of the eigendecomposition for dual-number
  matrices (a fact proved in the reference) and observing that the
  diagonal form of the block matrix
  \(\begin{pmatrix}0 & M^* \\ M & 0\end{pmatrix}\) is
  \(\Sigma \oplus (-\Sigma)\) where \(\Sigma\) is the normal form of
  \(M\) under SVD. This implies that \(\Sigma\) is unique.
\item
  See \cite{dualsvd}. Uniqueness can be reduced to
  the case of the involution \(z^* := z\).
\item
  This follows from \cite{contragredient}.
\item
  This follows from Theorem 7.2 of \cite{quaternions}.
\item
  See \cite{dualcomplex}.
\end{enumerate}

\textbf{Theorem 4}: \(\hermbar(R,*)\) admits an
isomorphism to a free abelian monoid with the following generators when
\((R,*)\) is any of the following \(*\)-rings:

\begin{longtable}{|l|l|l|l|}
  \hline
 & \(R\) & \(*:R \to R\) & Generators of \(\hermbar(R,*)\)\\
 \hline \hline
1 & \(0\) & \(z^* := z\) & \(\{[0]\}\) \\ \hline
2 & \(\mathbb R\) & \(z^* := z\) &
\(\{[x] : x \in \mathbb R; x > 0\} \cup \{[0]\}\) \\ \hline
3 & \(\mathbb C\) & \((a + bi)^* := a - bi\) &
\(\{[x] : x \in \mathbb R; x > 0\} \cup \{[0]\}\) \\ \hline
4 & \(\mathbb R[\varepsilon] / (\varepsilon^2)\) & \(z^* := z\) &
\(\begin{aligned}&\{[x + y\varepsilon] : x \in \mathbb R, y \in \mathbb R; x > 0\} \\&\cup \{[y \varepsilon] : y \in \mathbb R; y > 0\}\\&\cup \{[0]\}\end{aligned}\) \\ \hline
5 & \(\mathbb R[\varepsilon] / (\varepsilon^2)\) &
\((a + b\varepsilon)^* := a - b\varepsilon\) &
\(\begin{aligned}&\{[x] : x \in \mathbb R; x > 0\} \\&\cup \left\{\begin{bmatrix} x & -y\varepsilon \\ y\varepsilon & x\end{bmatrix} : x \in \mathbb R, y \in \mathbb R; x > 0, y > 0\right\}\\& \cup \{y\varepsilon : y \in \mathbb R; y > 0\}\\&\cup \{[0]\}\end{aligned}\) \\ \hline
6 & \(\mathbb C \oplus \mathbb C\) & \(\operatorname{swap}\) &
\(\{(J_m(re^{i\theta}),J_m(re^{i\theta})^T) : r \geq 0, \theta \in [0,\pi)\}\) \\ \hline
7 & \(\mathbb H\) & \(\begin{aligned}&(a + bi + cj + dk)^* :=\\& a - bi - cj - dk\end{aligned}\) &
\(\{[x] : x \in \mathbb R; x > 0\} \cup \{[0]\}\) \\ \hline
8 & \(\begin{aligned}&\operatorname{span}(1, i, \varepsilon j, \varepsilon k)\\&\subset \mathbb H[\varepsilon]/(\varepsilon^2)\end{aligned}\) &
\(\begin{aligned}&(a + bi + c\varepsilon j + d \varepsilon k)^* :=\\& a - bi - c\varepsilon j - d \varepsilon k\end{aligned}\)
&
\(\begin{aligned}&\{[x] : x \in \mathbb R; x > 0\} \\&\cup \left\{\begin{bmatrix} x & -\delta \\ \delta & x\end{bmatrix} : x \in \mathbb R; x > 0, \delta^2 = 0\right\}\\&\cup \{[0]\}\end{aligned}\) \\ \hline
9 & \(\mathbb Z\) & \(z^* := z\) & \begin{tabular}[t]{@{}l@{}}The generators are related to the set of \\adjacency matrices of all connected graphs. \\We don't give an explicit description.\end{tabular}\\ \hline
10 & \(\mathbb C\) & \(z^* := z\) & Very complicated. Given in \cite{complexsymmetric}. \\ \hline
\end{longtable}
\emph{Proof}. For each example in turn:

\begin{enumerate}
\item
  The ring here is the \(0\) ring, whose set of elements is \(\{0\}\).
  The definition of the operations \(\{+,-,\times\}\) is immediate.
\item
  A proof of this fact can be found in any undergraduate textbook on
  linear algebra for mathematics students, as long as that textbook
  covers the spectral theorem and SVD.
\item
  A proof of this fact can be found in any undergraduate textbook on
  linear algebra for mathematics students, as long as that textbook
  covers the spectral theorem and SVD.
\item
  See \cite{dualsvd}.
\item
  See \cite{dualsvd}.
\item
  Consequence of Theorem 2.
\item
  This follows from Theorem 7.2 of \cite{quaternions}.
\item
  See \cite{dualcomplex}.
\item
  In this case, we haven't stated the result clearly enough that we can
  prove it. If we wished to, we \emph{could} state it rigorously and
  prove it. Note that a unitary matrix \(U\) that is also an integer
  matrix is precisely a signed permutation matrix. Two undirected
  graphs, represented as adjacency matrix \(M\) and \(K\), are
  isomorphic iff there exists a permutation matrix \(P\) such that
  \(M = PKP^{-1}\). The generators are therefore connected graphs which
  are combined by direct sum to form unconnected graphs.
\item
  See \cite{complexsymmetric}.
\end{enumerate}

\section{Meta-conjectures}\label{meta-conjectures}

By a meta-conjecture, we mean a conjecture which is true for a large
class of rings, but not necessarily for all rings. We suspect that the
class of rings for which these conjectures are true is large.

\textbf{Conjecture 1}: For any ring \(*\)-ring \((R,*)\),
\(\hermbar(R,*)\) is free.

\textbf{Conjecture 2}: For any ring \(*\)-ring \((R,*)\),
\(\matbar(R,*)\) is free.

We know from theorem 1 that \(\sqbar(R)\) is free whenever
\(R\) is Artinian. \(\sqbar(R)\) may be free for some
non-Artinian rings as well. Therefore we do not pose this as one of the
meta-conjectures.

A proof of those two meta-conjectures may end up being non-constructive,
in the sense that it might fail to give the generators of the
corresponding free abelian monoids explicitly. We consider the search
for a constructive proof in some special cases to be worthwhile.

\bibliography{svds_via_monoids}{}
\bibliographystyle{plain}
\end{document}